\documentclass[11pt]{article}
\catcode`\@=11 \catcode`\@=12

\usepackage{amsmath}
\usepackage{amssymb}
\usepackage{amsthm}
\usepackage{oldgerm}

\newcommand{\Rm}{\mathbb{R}}
\newcommand{\RR}{\mathbb{R}}

\newcommand{\mT}{\ensuremath{\mathcal{T}}}

\newcommand{\mE}{\ensuremath{\mathcal{E}}}

\newcommand{\Qm}{\ensuremath{\mathbb{Q}}}
\newcommand{\Nm}{\ensuremath{\mathbb{N}}}

\newcommand{\mD}{\ensuremath{\mathcal{D}}}

\newcommand{\mO}{\ensuremath{\mathcal{O}}}

\newcommand{\vs}{\vspace{.2cm}}

\newtheorem{lem}{Lemma}
\newtheorem{thm}{Theorem}

\newtheorem{cor}[lem]{Corollary}
\newtheorem{prop}[lem]{Proposition}
\newtheorem{defn}[lem]{Definition}

\def\proof {\noindent{\sc{Proof. }}}
\def\qed {\mbox{}\hfill {\small \fbox{}} \\}
\def\lto{\longrightarrow}
\def\lmto{\longmapsto}

\def\leq{\leqslant}
\def\geq{\geqslant}

\title{The Monge problem for supercritical Ma\~n\'e potentials
on compact manifolds}
\author{Patrick  Bernard and Boris  Buffoni}
\date{ october 2005 }
\begin{document}

\maketitle
\begin{small}
Abstract:\\
We prove the existence of optimal transport maps for the Monge problem when the cost is a Finsler distance on a compact manifold. Our point of view consists in considering the distance as a Ma\~n\'e potential, and to rely on recent developments in the theory of viscosity solutions of the Hamilton-Jacobi equation.\\

R\'esum\'e:\\
On montre l'existence d'une application de transport optimale pour le probl\`eme de Monge lorsque le cout est une distance Finslerienne sur une vari\'et\'e compacte. Le nouveau point de vue consiste \`a consid\'erer la distance comme un potentiel de Ma\~n\'e, et \`a exploiter des d\'eveloppements r\'ecents sur les solutions de viscostit\'e de l'\'equation de Hamilton-Jacobi. \\

\vspace{1cm}

 Patrick Bernard\\
Institut Fourier, Grenoble,\\
on move to \\
CEREMADE\\
Universit\'e de Paris Dauphine\\
Pl. du Mar\'echal de Lattre de Tassigny\\
75775 Paris Cedex 16\\
France\\
\texttt{patrick.bernard@ceremade.dauphine.fr}\\

\vspace{1cm}

{Boris Buffoni}\\
School of Mathematics\\
\'Ecole Polytechnique F\'ed\'erale-Lausanne\\
SB/IACS/ANA Station 8\\
1015 Lausanne\\
Switzerland\\
\texttt{boris.buffoni@epfl.ch}

\end{small}

\newpage

The Monge transportation problem is to move one
distribution of mass into another in an optimal way. Before we
discuss this problem, let us describe a precise setting. We fix a
space $M$, which in the present paper will be a manifold, and a
cost function $c\in C(M\times M,\Rm)$. Given two probability
measures $\mu_0$ and $\mu_1$, we call  transport map a Borel map
$F:M\lto M$ that transports $\mu_0$ onto $\mu_1$. An optimal
transport map is a transport map $F$ that minimizes the total
cost
$$
\int_M c(x,F(x))d\mu_0
$$
among all transport maps. In many situations,  optimal transport
maps have remarkable  geometric properties, at least at a formal
level. Some of these properties were investigated by Monge at the
end of the 
eighteenth century.

The question 
of existence of optimal transport maps was
discussed much later in the literature. Some major steps were
made by Kantorovich in 1942. He introduced both a relaxed
problem and a dual problem that opened new approaches to the
existence problem. When the cost is the square of the distance on
an Euclidean vector space, Brenier proved the existence of an
optimal transport map  in \cite{Br:91} and also provided an
interesting geometric description  on the optimal maps, which have 
to be the gradient of a convex function. The argument was
simplified, taking advantage of the Kantorovich dual problem, by
Gangbo, \cite{Ga:94}, and  extended in many directions by Gangbo,
McCann, \cite{GaMc:96}  and other authors, see our paper
\cite{BeBu:first} for more details.

The case where the cost function is the distance on an Euclidean
vector space is very natural, but more difficult. Sudakov
announced a proof of the existence of an optimal map in 1979, but
a gap was recently found in this proof. The strategy was to
decompose the space into pieces of smaller dimension on which
transport maps can be more easily built, and to glue these maps
together. An essential hypothesis of the result is that the
measure $\mu_0$ is absolutely continuous. In the construction, it
is necessary to control how this hypothesis behaves under
decomposition to the subsets. Sudakov was not aware of these
difficulties, and made wrong statements  at that point, as was
discovered only much later. It is interesting to notice for
comparison  that similar kind of difficulties had been faced and
solved ten years earlier by Anosov in his ergodic theory of
hyperbolic diffeomorphisms.

Correct proofs in the spirit of the work of Sudakov were written
simultaneously by Caffarelli, Feldman and McCann in
\cite{CaFeMc:02},  Trudinger and Wang in \cite{TrWa:01} and
slightly afterwards by Ambrosio, \cite{Am:00}.
   These authors manage to build  a decomposition
of the space into line segments that have to be preserved by
transport maps, called transport rays. They prove that the
direction of these rays vary Lipschitz continuously. This
regularity implies that the absolutely continuous measure $\mu_0$
has absolutely continuous decompositions on these rays. See
\cite{Am:00} for a remarkably written discussion  on these works
and of Sudakov's mistake. Before the proof of Sudakov were
completed in these papers, Evans and Gangbo had provided a
different proof under more stringent hypotheses in \cite{EvGa:99}.
This proof is long and complicated, but it now appears as
 the first proof of the existence of a transport map
in the case where the cost is a  distance.

The methods inspired from Sudakov seem to allow many kinds of
generalizations. The paper  \cite{CaFeMc:02} treats all norms
whose unit ball is smooth and strictly convex. It is worth
mentioning that flat part in the unit ball represent a major
difficulty. An important progress  have been recently made in
\cite{AKP:04}, which studies norms whose unit ball is a
polyhedron. In another direction,  Feldman and McCann
\cite{FeMc:02} have treated the case where the cost is the
distance on a Riemannian manifold. In this generalized setting,
transport rays are not any more  line segments, but pieces of
geodesics. It is in this direction we will pursue in the present
text.

Our goal is to prove the existence of
 transport maps for Finsler distances on 
manifolds
(possibly non symmetric distances). In order to avoid superficial 
(and less superficial)  additional
technicalities, we shall work on a compact manifold $M$. Our
first  novelty is a new approach of the geometric part of the
proof, that is the decomposition into transport rays. We believe
that this new approach is interesting because, beyond  being more
general, it enlightens new links between the Monge problem and the
general theory of Hamilton-Jacobi equations as presented by
Fathi  in \cite{Fa:un}. In fact, all the relevant properties of
the decomposition into transport rays are obtained by
straightforward applications of results of \cite{Fa:un}. In order
to finish the proof, we rely on a secondary
 variational principle, in the lines of \cite{AmPr:03} and
\cite{AKP:04}. Our treatment of this secondary principle is quite
different from these papers, and it is, we believe,
shorter and clearer than the methods previously used in the literature.

This paper was born during the visit of the first author to the
Bernoulli center in EPFL, Lausanne, in Summer 2003. We wish to
thank this institution for its support.

\section{Introduction}
We state two versions of our main result,
and prove the equivalence between these two statements.

\subsection{Optimal transport maps for Finsler distances}
In  the present paper, the space $M$ is a smooth compact connected  manifold
without boundary. 
We equip the manifold $M$ with a $C^2$ Finsler metric, that is 
 the data,
for each $x\in M$, of a non-negative convex
function $v\lmto \|v\|_x$ on $T_xM$  such that
\begin{itemize}
\item $\|\lambda v\|_x=\lambda\|v\|_x$ for all $\lambda>0$ and all $v\in T_xM$
(positive 1-homogeneity),
\item
the function $(x,v)\lmto \|v\|_x$ is $C^2$ outside of the zero section, 
\item
for each $x\in M$, the function $v\lmto \|v\|^2_x$ has positive definite
Hessian at all vectors $v\neq 0$.
\end{itemize}
As a consequence, $\|v\|_x=0$ exactly when $v=0$ (positivity).
Note that $\|\cdot\|_x$ is not assumed symmetric
 (that is, $\|-v\|_x\neq \|v\|_x$ is allowed).
In standard terminology, the function  $v\lmto \|v\|_x$
is a Minkowsky metric on the vector space $T_xM$.
See \cite{BCS} for more one Finsler metrics
(in particular Theorem 1.2.2 and paragraph 6.2).
We define the 
length of each smooth curve $\gamma:[0,T]\lto M$
by the expression
$$
l(\gamma)=\int_0^T\|\dot \gamma(t)\|_{\gamma(t)}dt.
$$
The Finsler distance $c$ is then given by the expression
$$
c(x,y)=\inf _{\gamma} l(\gamma)
$$
where the infimum is taken on the set of smooth curves 
$\gamma:[0,T]\lto M$ 
(where $T$ is any positive number) 
which satisfy $\gamma(0)=x$ and $\gamma(T)=y$.
Note that the value of the 
 infimum would not be changed  by  imposing the additional
 requirement that 
   $\|\dot \gamma(t)\|_{\gamma(t)}\equiv 1$.
The Finsler distance $c(x,y)$ is not 
necessarily symmetric,
and thus is not properly speaking a distance. It does satisfy
the triangle inequality, and $c(x,y)=0$ if and only if
$x=y$. 

We shall consider the Monge transportation problem for the cost
$c$. Given a Borel measure $\mu_0$ on $M$, and a Borel map
$F:M\lto M$, we define the image measure $F_{\sharp}\mu_0$ by
$$
F_{\sharp}\mu_0(A):= \mu_0(F^{-1}(A))
$$
for each Borel set $A\subset  M$. The map $F$ is said to
transport $\mu_0$ onto $\mu_1$ if $F_{\sharp}\mu_0=\mu_1$. 
We will present a short proof of:

\begin{thm}
Let $c(x,y)$ be a Finsler distance on $M$.
Let $\mu_0$ and $\mu_1$ be two
probability measures on $M$, such that $\mu_0$ is absolutely
continuous with respect to the Lebesgue class. Then there exists
a Borel map $F:M\lto M$ such that
$F_{\sharp}\mu_0=\mu_1$, and such that the inequality
$$
\int_M c(x,F(x))d\mu_0 \leq \int_M c(x,G(x))d\mu_0
$$
holds for each Borel map  $G:M\lto M$ satisfying
$G_{\sharp}\mu_0=\mu_1$. In other words, there exists an optimal
map for the Monge transportation problem.
\end{thm}

\subsection{Optimal transport maps for supercritical Ma\~n\'e potentials}
We shall now present a generalization of Theorem 1,
which is the natural setting for our proof.
A function $L:TM\lto \Rm$ is called a Tonelli Lagrangian
if it is $C^2$ and satisfies:\\
\textbf{Convexity} For each $x\in M$, the function $v\lmto
L(x,v)$ is convex with positive definite Hessian
at each point.\\
\textbf{Superlinearity} For each  $x\in M$, we have
$L(x,v)/\|v\|_x\lto \infty$
as $\|v\|_x\lto \infty$.\vs\\
Given a Tonelli Lagrangian $L$ and a time 
$T\in ]0,\infty)$, we define the cost function
$$
c^L_T(x,y)=\min_{\gamma} \int _0^T L(\gamma(t), \dot \gamma(t)) dt
$$
where the minimum is taken on the set of curves $\gamma\in
C^2([0,T],M)$ satisfying $\gamma(0)=x$ and $\gamma(T)=y$. That
this minimum exists is  standard,
see \cite{Ma:91} or \cite{Fa:un}.
The function 
$$c^L(x,y):=\inf _{T\in ]0,\infty)} c_T(x,y)
$$
is called the Ma\~n\'e potential of $L$.
It was introduced and studied by Ricardo Ma\~n\'e and then his
students in \cite{Ma:97,CoDeIt}. Without additional hypothesis,
the Ma\~n\'e potential may be identically $-\infty$. So we assume
in addition:
\\
\textbf{Supercriticality}
For each $x \neq y\in M^2$, we have $c(x,y)+c(y,x)>0$.\vs\\
 The following result of Ma\~n\'e
\cite{Ma:97,CoDeIt}  makes this hypothesis natural.

\begin{prop}
Let $L\in C^2(TM,\Rm)$ be a Tonelli Lagrangian.
For $k\in \Rm$, let $c^{L+k}$ be the Ma\~n\'e potential associated to
the Lagrangian $L+k$. There exists a constant $k_0$ such that
\begin{itemize}
\item For $k< k_0$, then  $c^{L+k}\equiv -\infty$ and  the Lagrangian $L+k$
is called subcritical.
\item For $k\geq k_0$, the Ma\~n\'e potential
 $c^{L+k}$ is a Lipschitz  function
on $M\times M$  that satisfies the triangle inequality
$$
c^{L+k}(x,z)\leq c^{L+k}(x,y)+c^{L+k}(y,z)
$$
for all $x,y$ and $z$ in $M$. In addition, we have $c^{L+k}(x,x)=0$
for all $x\in M$.
\item
For $k>k_0$, the Lagrangian $L$ is supercritical, which means that
$c^{L+k}(x,y)+c^{L+k}(y,x)>0$ for $x\neq y$ in $M$.
\end{itemize}
\end{prop}
We will explain that the following theorem 
is equivalent to Theorem 1.

\begin{thm}
Let $c^L(x,y)$ be the Ma\~n\'e potential associated to a
supercritical Tonelli Lagrangian $L$. Let $\mu_0$ and $\mu_1$ be two
probability measures on $M$, such that $\mu_0$ is absolutely
continuous with respect to the Lebesgue class.
Then there exists an optimal transport map for the Monge transportation problem
with cost  $c^L(x,y)$.
\end{thm}

\subsection{Supercritical Ma\~n\'e potentials and Finsler distances}
We prove the equivalence between Theorem 1 and Theorem 2.
%
%
%
To each Tonelli Lagrangian $L$, we associate  the Hamiltonian 
 $H\in C^2(T^*M,\Rm)$ defined by 
$$
H(x,p)=\max_{v\in T_xM} p(v)-L(x,v)
$$
and the energy function 
$E\in C^2(TM,\Rm)$ defined by
$$
E(x,v)=\partial_vL(x,v).v-L(x,v).
$$
The  function $H$ is also convex and superlinear.
The mapping
$\partial_v L:TM\lto T^*M$ is a $C^1$ diffeomorphism, whose
inverse is the mapping $\partial_p H$.
We have $E=H\circ \partial_v L$.

\begin{lem}\label{energy}
Let $L$ be a supercritical Tonelli  Lagrangian. There exists a constant
$K$ such that, for each $x\neq y$ in $M$, there exists a time
$T\in ]0,K]$ and a minimizing extremal $\gamma\in C^2([0,T],M)$
such that $\gamma(0)=x$, $\gamma(T)=y$, and 
$\int_0^TL(\gamma(t)
,\dot \gamma(t))dt =c_T^L(x,y)=c^L(x,y)$.
Moreover, if $\gamma$ is such a curve, then 
$$
E(\gamma(t),\dot\gamma(t))\equiv 0.
$$
\end{lem}
\proof
We shall prove, and use, this lemma only in the case where $L$ is
positive.
Note first that the function $(x,y,T)\lmto c_T^L(x,y)$
is continuous on $M\times M\times ]0,\infty)$.
It is not hard to see, in view of the superlinearity of $L$, that
the function $T\lmto c^L_T(x,y)$ goes to infinity
as $T$ goes to zero if $x\neq y$.
On the other hand, setting $\delta=\inf L>0$, we obviously have the minoration
$c^L_T\geq \delta T$.
Since $c^L_1$ is bounded, this implies the existence of a constant $K$
such that $c^L_T> c^L_1$ for $T\geq K$.
As a consequence, the function $T\lmto c_T(x,y)$
reaches its minimum on $]0,K]$ for each $x\neq y$.

Let $x\neq y$ be two points on $M$.
There exists a $T\in ]0,K]$
such that $c^L_T(x,y)=c^L(x,y)$.
Now by standard results on the calculus of variations, 
there exists a $C^2$  curve 
 $\gamma:[0,T]\lto M$  satisfying $\gamma(0)=x$,
$\gamma(T)=y$ and $\int_0^TL(\gamma(t),\dot \gamma(t))dt =c^L(x,y)$.
In addition, this curve satisfies the Euler-Lagrange equations, 
and in particular the energy $E(\gamma(t),\dot \gamma(t))$
is constant on $[0,T]$.

Let $\gamma_{\lambda}:[0,\lambda T]\lto M$
be defined by $\gamma_{\lambda}(t)=\gamma(t/\lambda)$.
The function 
$$f(\lambda):=
\int_0^{T\lambda}L(\gamma_{\lambda},\dot \gamma_{\lambda})dt
=\lambda\int_0^{T}L(\gamma,\lambda^{-1}\dot \gamma)dt$$
clearly has to reach its minimum at $\lambda =1$.
On the other hand, a classical computation shows that 
the function $f$ is differentiable, and that 
$f'(1)=-\int_0^T E(\gamma(t),\dot \gamma(t))dt$.
This proves that $E(\gamma(t),\dot \gamma(t))\equiv 0$.
\qed

The following proposition implies that the transportation
problem for Finsler distances, the transportation problem
for supercritical Ma\~n\'e potentials, and the transportation problem
for the Ma\~n\'e potentials of positive Tonelli Lagrangians
are equivalent problems.

\begin{prop}\label{costs}
If $L$ is a supercritical Tonelli Lagrangian,
then there exists  a Finsler distance $c$ 
(associated to a $C^2$ Finsler metric)
and a smooth function $f:M\lto \Rm$ such that 
$$
 c(x,y) =c^L(x,y)+f(y)-f(x).
$$
Conversely, given a Finsler distance $c$ 
(associated to a $C^2$ Finsler metric)
there exists a positive Tonelli Lagrangian $L$
such that 
$$
 c(x,y) =c^L(x,y).
$$
\end{prop}
\proof
The first part of this proposition
is the content of  \cite{ISM:00}.
For the converse, we consider the Lagrangian
$$\tilde L(x,v)=
\frac{1+\|v\|_x^2}{2}.
$$
Note that the associated energy function
is 
$$
\tilde E(x,v)=\frac{\|v\|_x^2-1}{2}.
$$

Let us now consider a positive Tonelli Lagrangian $L$
such that 
$0<L\leq \tilde L$
and such that $L=\tilde L$ on the set $\{\|v\|_x\geq 1/2\}\subset TM$.
In order to see that such a Lagrangian exists,
consider a smooth convex function $f:[0,\infty)\rightarrow [0,\infty)$
that vanishes in a small neighborhood of $0$ and such that $f(s)\leq (1+s^2)/2$
with equality for $s\geq 1/2$;
and consider a smooth function $g(x,v): TM\lto \Rm$
such that $g$ is positive and $\partial_2g(x,v)$ is positive
definite  where $f(\|v\|_x)$ vanishes,
and such that $g$ is zero on $\{\|v\|_x\geq 1/2\}\subset TM$.
It is easy to check that the Lagrangian 
 $L(x,v)=f(\|v\|_x)+\epsilon g(x,v)$ satisfies
 the desired requirements  when $\epsilon>0$
is small enough,

Since $L\leq \tilde L$, we have $c^L\leq c^{\tilde L}$.
Moreover, for $x\neq y$, let $\gamma:[0,T]\lto M$
be the optimal (for $c^L$) trajectory obtained by lemma \ref{energy}.
We have $E(\gamma(t),\dot \gamma(t))=
\tilde E(\gamma(t),\dot \gamma(t))\equiv 0$
hence $\|\dot \gamma\|\equiv 1$. As a consequence, 
$$
c^L(x,y)=\int_0^TL(\gamma(t),\dot \gamma(t))dt=
\int_0^T \tilde L(\gamma(t),\dot \gamma(t))dt
=l(\gamma)\geq c(x,y).
$$

We have proved that 
$c\leq c^L\leq c^{\tilde L}$. These are equalities because $c^{\tilde L}\leq c$. Indeed, for all
smooth curves
$\gamma:[0,T]\lto M$ which satisfy $\gamma(0)=x$, $\gamma(T)=y$ and
$\|\dot \gamma(t)\|_{\gamma(t)}\equiv 1$ with $T=l(\gamma)>0$, we get
$c^{\widetilde L}(x,y)
\leq \int_0^T \tilde L(\gamma(t),\dot \gamma(t))dt
= l(\gamma)$.
Since $c(x,y)$ is the infimum of the lengths of such curves, we have 
$c^{\tilde L}(x,y)\leq c(x,y)$.
\qed

\subsection{General convention}\label{convention}
In the sequel, we will prove Theorem 2 for
a positive Lagrangian $L$ and 
denote $c^L$ simply by $c$. 
In view of Proposition \ref{costs},
this implies the general form of Theorem 2, as well as 
Theorem 1.
We fix, once and for all, a positive Tonelli
Lagrangian $L$, and a positive number $\delta>0$
such that
$$ L(x,v)\geq \delta$$
for each $(x,v)\in TM$.

The general scheme of our  proof is somewhat similar to the one
introduced by Sudakov, and followed, in \cite{CaFeMc:02},
\cite{TrWa:01}, \cite{Am:00}, \cite{FeMc:02}, \cite{AKP:04} and
other papers. Like these papers, our proof involves a
decomposition of the space $M$ into distinguished curves, called
transport rays. We introduce these rays in section \ref{rays} and
describe their geometric properties. In this geometric part of
the proof, our point of view is quite different from the
literature as we emphasize the link with  the theory of viscosity
solutions as developed in \cite{Fa:un}, and manage to obtain all
the relevant properties of transport rays as a straightforward
application of general results of \cite{Fa:un}. For the second
part of the proof, all the papers mentioned above involve subtle
decompositions of measures on these transport rays. It is at this
step that the paper of Sudakov contains a gap. We simplify this
step by introducing a secondary variational principle in section
\ref{plan}. Note that secondary variational principles have
already been introduced by Ambrosio, Kirchheim  and Pratelli in
\cite{AKP:04} for related problems. This secondary problem is
studied in section \ref{proof} by a quite simple method, which,
surprisingly, seems new. This methods allows a neat clarification
of the end of the proof compared to the existing literature. All
the difficulties involving  measurability issues  and absolute
continuity of disintegrated measures are reduced to a single and
simple Fubini-like theorem, exposed in section \ref{fubini}.

\section{Transport plans}\label{plan}
We introduce our secondary variational principle and recall the
necessary generalities on the Monge problem. Beyond the
references provided below, the
 pedagogical texts \cite{Am:00,RaRu:98,Vi:03} may help the reader
who wants more details. We  define the quantity
$$
C(\mu_0,\mu_1):=\inf_F \int _M c(x,F(x)) d\mu_0,
$$
where the infimum is taken on the set of Borel maps $F:M\lto M$
that transport $\mu_0$ onto $\mu_1$. It is useful, following
Kantorovich,
 to relax this infimum to a nicer minimization problem.
A Borel measure $\mu$ on $M\times M$ is called a transport plan
if it satisfies the equalities $\pi_{i\sharp}\mu=\mu_i$, where
$\pi_0:M\times M\lto M$ is the projection on the first factor, and
$\pi_1:M\times M\lto M$ is the projection on the second factor.
Clearly, any transport map $F$ can be considered as the transport
plan
$$(Id\times F)_{\sharp}\mu_0.$$
Following Kantorovich, we consider the minimum
$$
K(\mu_0,\mu_1)=\min_{\mu} \int_{M\times M} cd\mu
$$
taken on the set of transport plans. It is well-known and easy to
prove that this minimum exists. The equality
$$
K(\mu_0,\mu_1)=C(\mu_0,\mu_1)
$$
holds if $\mu_0$ has no atom, see \cite{Am:00}, Theorem 2.1. Note
that this equality is very general, see \cite{Pr:gen} for a
discussion. As a consequence, if $\mu_0$ has no atom, it is
equivalent to prove the existence of an optimal transport map and
to prove that there exists an optimal transport plan
concentrated  on the graph of a Borel function.

Let us define the second cost  function
$$
\sigma(x,y)=(c(x,y))^2.
$$
This cost is chosen in order that the following refined form
of Theorem 1 holds.

\begin{thm}
Let $\mO$ be the set of optimal transport plans for
$K(\mu_0,\mu_1)$ with the cost $c$. The minimum
$$\min _{\mu\in \mO}\int \sigma d\mu$$
exists. In addition, if $\mu_0$ is absolutely continuous, then
there is one and only one transport plan $\mu$ realizing this
optimum, and this transport plan is concentrated on the graph of
a Borel function that is an optimal transport map for the cost
$c$.
\end{thm}

This result will be proved in section \ref{proof}. The idea of
introducing secondary variational problem as in this statement
has already been used by Ambrosio,  Kirchheim,  and  Pratelli, see
\cite{AKP:04} and also \cite{AmPr:03}. Our treatment  in section
\ref{proof} is inspired from these references,
although it is somewhat different. It allows
substantial simplifications compared to the literature.

\section{Kantorovich potential and calibrated curves}\label{rays}
We present the decomposition in transport ray, which is the
standard initial step in the construction of optimal maps. This
construction is based on well-understood general results on
viscosity sub-solutions of the Hamilton-Jacobi equation, as
presented in \cite{Fa:un}. Making this connection is one of the
novelties of the present paper.

Since  the cost function we consider satisfies the triangle
inequality
$$
c(x,z)\leq c(x,y)+c(y,z)
$$
for all $x,y$ and $z$ in $M$, as well as the identity $c(x,x)=0$
for all $x\in M$, we can take advantage of  the following general
duality result, inspired from Kantorovich, see for example
\cite{Vi:03}, \cite{Ev:99} and  \cite{FeMc:02}.

\begin{prop}\label{Kanto}
Given two measures $\mu_0$ and $\mu_1$, there exists a function
$u\in C(M,\Rm)$ that satisfies
$$
u(y)-u(x)\leq c(x,y)
$$
for all $x$ and $y$ in $M$, and
$$
K(\mu_0,\mu_1)=\int_M ud(\mu_1-\mu_0).
$$
In addition, for each optimal transport plan $\mu$, the equality
$u(y)-u(x)= c(x,y)$ holds for $\mu$-almost every $(x,y)\in M^2$.
The function $u$ is called a Kantorovich potential.
\end{prop}

The present paper is born from the observation that the
Kantorovich potentials are viscosity subsolutions of the
Hamilton-Jacobi equation as studied in \cite{Fa:un}.
In order to explain this connection, it is necessary to define
the Hamiltonian function 
 $H\in C^2(T^*M,\Rm)$ by
$$
H(x,p)=\max_{v\in T_xM} p(v)-L(x,v).
$$
Note that the mapping
$\partial_v L:TM\lto T^*M$ is a $C^1$ diffeomorphism, whose
inverse is the mapping $\partial_p H$. 
It is proved in \cite{Fa:un} that the following
properties are equivalent for a function $w\in C(M,\Rm)$.

\begin{enumerate}
\item
The function $w$ satisfies the inequality $w(y)-w(x)\leq c(x,y)$
for all $x$ and $y$ in $M$.
\item
The function $w$ is a viscosity sub-solution of the
Hamilton-Jacobi equation $H(x,dw)=0$, \textit{i.e.}  each smooth
function $f:M\lto \Rm$ satisfies the inequality $H(x,df(x))\leq 0$
at each point of minimum $x$ of the difference $f-w$.
\item
The function $w$ is Lipschitz and satisfies the inequation
$H(x,dw_x)\leq0$  at almost every point.
This inequality then holds at all point of differentiability $x$ of $w$.
\item
The function $w$ is Lipschitz and, for almost every $x\in M$, it
 satisfies the inequation
$\forall v\in T_xM~L(x,v)\geq dw_x(v)$.
This inequality then holds at all point of differentiability $x$ of $w$.
\end{enumerate}

\begin{itshape}
Although there may exist several Kantorovich potentials, we
shall fix one of them, $u$, for the sequel.
\end{itshape}

\begin{defn}
{\em
Following Fathi  \cite{Fa:un}, we call \textit{calibrated curve}
a continuous 
and piecewise $C^1$ curve 
$\gamma:I\lto M$ that satisfies

\begin{equation}\label{calibrated}
u(\gamma(t))-u(\gamma(s))=\int_s^t L(\gamma(\tau),\dot \gamma(\tau))d\tau
=c(\gamma(s),\gamma(t))
\end{equation}
whenever $s\leq t$ in $I$, where $I$ is a non empty interval of $\Rm$
(possibly a point).
A calibrated curve
$\gamma:I\lto M$ is called {\em non-trivial} if the interval $I$ has
non-empty interior. 
}\end{defn}
Note that the first of the equalities in \eqref{calibrated} 
implies the second, since
the inequalities

$$
u(\gamma(t))-u(\gamma(s))\leq c(\gamma(s),\gamma(t)) \leq
\int_s^t L(\gamma(\tau),\dot \gamma(\tau))d\tau
$$
hold for any curve $\gamma$. It is obvious that non-trivial calibrated curves
are minimizing extremals of $L$, and as a consequence they are
$C^2$ curves. In addition, the concatenation of two calibrated
curves is calibrated, so that each calibrated curve can be
extended to a maximal calibrated curve, that is,  its domain $I$ 
cannot be further extended without loosing calibration. Note
that $I$ is closed when $\gamma$ is a maximal calibrated curve. 

\begin{defn}
{\em
We call \textit{transport ray} the image of a
non-trivial maximal calibrated curve.
}\end{defn}

It will be useful to consider, following \cite{CaFeMc:02} and
\cite{FeMc:02}, the functions $\alpha$ and $\beta:M\lto
[0,\infty)$ defined as follows:
\begin{itemize}
\item
$\alpha(x)$ is the supremum of all times $T\geq0$ such that there
exists a calibrated curve $\gamma:[-T,0]\lto M$ that satisfies
$\gamma(0)=x$.
\item
$\beta(x)$ is the supremum of all times $T\geq0$ such that there
exists a calibrated curve $\gamma:[0,T]\lto M$ that satisfies
$\gamma(0)=x$ 
(the fact that $\alpha$ and $\beta$ are finite is a consequence of
Lemma~\ref{increase} below).
\end{itemize}

\begin{defn}{\em
Let us denote by $\mT$ the subset of $M$ obtained as the union of
all transport rays, 
or equivalently the set of points $x\in M$
such that $\alpha(x)+\beta(x)>0$. For
$\epsilon \geq 0$, we denote by $\mT_{\epsilon}$ the set of
points $x\in M$ that satisfy $\alpha(x)>\epsilon$ and
$\beta(x)>\epsilon$.  Clearly $\mT_{\epsilon}\subset \mT$ for all 
$\epsilon\geq 0$.
The set $\mE:=\mT-\mT_0$ is the set of ray ends.
}\end{defn}

\begin{prop}\label{Lip}
The function $u$ is differentiable at each point of $\mT_0$. For
each point $x\in \mT_0$, there exists a single maximal calibrated
curve 
\begin{equation}\label{def: gammax}
\gamma_x:[-\alpha(x),\beta(x)]\lto M~~\hbox{such that}~~
\gamma_x(0)=x. 
\end{equation}
This curve satisfies the relations
$$
du_x =\partial_vL(x,\dot \gamma_x(0)) \text{ or equivalentely }
\dot \gamma_x(0)=\partial_pH(x,du_x).
$$
For each $\epsilon>0$, the differential $x\lmto du_x$ is
Lipschitz on $\mT_{\epsilon}$, or equivalently the map $x\lmto
\dot\gamma_x(0)$ is Lipschitz on $\mT_{\epsilon}$.
\end{prop}

\proof This proposition is Theorem 4.5.5 of Fathi's book
\cite{Fa:un}. \qed

\begin{lem}\label{increase}
Let $\gamma:[a,b]\lto M$ be a non-trivial calibrated curve. 
Then, for all $t\in ]a,b[$, the function $u$ is differentiable at 
$\gamma(t)$ and
$$du_{\gamma(t)}(\dot \gamma(t))\geq \delta,$$
(see \ref{convention} for the definition of $\delta$).
As a consequence, the map $\gamma:[a,b]\lto M$ is an embedding
(it is one to one and has non-zero derivative on $]a,b[$) 
and transport rays are non-trivial embedded arcs.
 \end{lem}

\proof 
Since $u$ is a viscosity sub-solution of the Hamilton-Jacobi equation,
see the equivalence below Proposition \ref{Kanto},
we have 
$L(x,v)\geq du_x(v)$ for all $v\in T_xM$
at each point of differentiability of $u$.
As a consequence, the inequality 
$$
L(\gamma(t),\dot \gamma(t))
\geq du_{\gamma(t)}(\dot \gamma(t))
$$
holds for each $t\in ]a,b[$.
Integrating the above
inequality gives
$$
\int_{a}^{b}L(\gamma(t),\dot \gamma(t))dt \geq
u(\gamma(b))-u(\gamma(a)),
$$
which is an equality because the curve $\gamma$ is calibrated. As
a consequence, we have
$$
du_{\gamma(t)}(\dot \gamma(t))=L(\gamma(t),\dot \gamma(t))
\geq \delta
$$
for all $t\in ]a,b[$.\qed

\begin{lem}
The functions $\alpha$ and $\beta$ are bounded  and upper
semi-continuous, hence Borel measurable.
As a consequence, the sets $\mT$ and $\mT_{\epsilon}$,
$\epsilon \geq 0$, are Borel.
\end{lem}
\proof We shall consider only the function $\alpha$. 
We have just seen
that, for each non-trivial calibrated curve $\gamma:I\lto M$, the function
$f(t)=u\circ \gamma(t)$ is differentiable and satisfies
$f'(t)\geq \delta$. Since the continuous function $u$ is
bounded on the compact manifold $M$, we conclude that the
functions $\alpha$ and $\beta$ are bounded. 
In order to prove
that the function $\alpha$ is semi-continuous, let us consider a
sequence $x_n\in M$ that is converging to a limit $x$ and is
such that $\alpha(x_n)\geq T$. We have to prove that
$\alpha(x)\geq T$. There exists a sequence $\gamma_n:[-T,0]\lto M$
of calibrated curves such that $\gamma_n(0)=x_n$. There exists a
subsequence of $\gamma_n$ that is converging uniformly on
$[-T,0]$ to a curve $\gamma:[-T,0]\lto M$. It is easy to see that
the curve $\gamma$ is calibrated and satisfies $\gamma(0)=x$. As
a consequence, we have $\alpha(x)\geq T$.
 \qed

\begin{defn}
{\em
For $x\in M$, let us denote by $R_x$ the union of the transport
rays containing $x$. 
We also denote by $R_x^+$ the set of points
$y\in M$ such that $u(y)-u(x)=c(x,y)$.
}\end{defn}
Note that $R_x=\gamma_x([-\alpha(x),\beta(x)])$
when $x\in \mT_0$, where $\gamma_x$ is given in \eqref{def: gammax}.
\begin{lem}
 We have $R_x^+=\gamma_x([0,\beta(x)])$
when $x\in \mT_0$, and $R^+_x=\{x\}$ when $x\in M-\mT$.
\end{lem}
\proof Let $x$ be a point of $\mT_0$. By the calibration property
of
 $\gamma_x$, we have, for $t\in [0,\beta(x)]$,
$\gamma_x(t)-\gamma_x(0)=c(\gamma_x(0),\gamma_x(t))$, which is
precisely saying that $\gamma_x(t)\in R_x^+$. Conversely, let us
fix a point $x\in M$ and let $y$ be a point of $R_x^+$. There
exists a time $T\geq 0$ and a curve $\gamma:[0,T]\lto M$ such that
$\int_0^T L(\gamma(t),\dot \gamma(t)) dt =c(x,y)$,
 $\gamma(0)=x$ and $\gamma(T)=y$.
Since $c(x,y)=u(y)-u(x)$, the curve $\gamma$ is calibrated. If
$x\in \mT_0$, then $\gamma=\gamma_{x}\big|_{[0,T]}$ hence
$y=\gamma(T)=\gamma_x(T)\in \gamma_x([0,\beta(x)])$. If $x\not
\in \mT$, then there is no nontrivial calibrated curve starting
at $x$, so we must have $y=x$ in the above discussion, and
$R_x^+=\{x\}$.
\qed

\begin{prop}\label{optimality}
The transport plan $\mu$ is optimal for the cost $c$ if and only
if it is concentrated on the closed set 
$$\cup _{x\in M}\{x\}\times R_x^+=\{(x,y)\in M^2: c(x,y)=u(y)-u(x)\}.$$
\end{prop}
\proof
Recall from Proposition \ref{Kanto} that any optimal transport
plan is concentrated on this set. Reciprocally, if $\mu$ is
a transport plan concentrated on this set, then
$$K(\mu_0,\mu_1)=\int_M ud(\mu_1-\mu_0)
=\int_{M^2}(u(y)-u(x))d\mu
=\int_{M^2}cd\mu$$
and $\mu$ is thus optimal.
\qed

\section{Fubini Theorem}\label{fubini}
The geometric informations on transport rays that have been obtained in the
preceding section imply the following crucial Fubini-like result, whose
proof is the goal of the present section :

\begin{prop}
\label{Fubini-like}
Let $\Lambda$ be a Borel  subset of $\mT$ such that the intersection
$\Lambda\cap R$ has zero $1$-Hausdorff measure for each transport
ray $R$. Then the set $\Lambda$ has zero Lebesgue measure.
\end{prop}
For comparison with the literature, see
\cite{CaFeMc:02} Lemma 25, \cite{TrWa:01} section 4 or
\cite{FeMc:02} Lemma 24,
 we mention:
\begin{cor}
The set $\mE=\mT-\mT_0$ of ray ends has zero Lebesgue measure.
\end{cor}
We now turn to the proof of Proposition \ref{Fubini-like},
which occupies the end of this section. 
The method is standard. Let $k$ be the dimension of $M$.

\begin{defn}
{\em
We call \textit{transport beam} 
any given pair $(B,\chi)$, 
where $B$ is a bounded Borel subset of $\Rm^k$
and $\chi:B\lto M$ is a Lipschitz map (not necessarily one-to-one)
such that:
\begin{itemize}
\item There exists a bounded Borel  set $\Omega\in \Rm^{k-1}$
and two 
bounded Borel functions $a<b:\Omega\lto \Rm$ such that 
$$
B=\{(\omega,s)\in \Omega \times \Rm \text{ s. t. }
a(\omega)\leq s\leq b(\omega)\}
\subset  \Rm^k= \Rm^{k-1}\times \Rm
$$
\item
For each $\omega \in \Omega$, the curve 
$\chi_{\omega}:[a(\omega),b(\omega)]\lto M$ given by
$\chi_{\omega}(s)=\chi(\omega,s)$
is a  calibrated curve.
\end{itemize}
}\end{defn}

\begin{lem}
If $(B,\chi)$
is a transport beam, then the set $\Lambda \cap \chi(B)$
has zero Lebesgue measure.
\end{lem}

\proof
For each $\omega \in \Omega$, the curve $\chi_{\omega}$
is a bilipschitz homeomorphism onto its image.
Since in addition, the set 
$\Lambda \cap \chi(\{\omega\}\times [a(\omega),b(\omega)])$
has zero $1$-Hausdorff measure,
the set  $\chi^{-1}(\Lambda)$ intersects each vertical 
line 
$\{\omega \}\times \Rm$ along a set of zero 1-Hausdorff measure.
In view of the classical Fubini theorem,
the set $\chi^{-1}(\Lambda)$ has zero Lebesgue
 measure in $\Rm^k$.
Since the $k$-Hausdorff measure on $B$
(associated to the restricted Euclidean metric)
is the restriction to $B$ of the Lebesgue measure of $\Rm^k$,
the set $\chi^{-1}(\Lambda)$
has  zero  $k$-Hausdorff measure in 
$B$, for the Hausdorff measure associated to the induced metric.
Since Lipschitz maps send sets of zero  $k$-Hausdorff
measure onto  sets of zero  $k$-Hausdorff
measure,
we conclude that the set 
$\Lambda \cap \chi(B)\subset 
\chi(\chi^{-1}(\Lambda)) $
has zero Lebesgue measure in $M$.
\qed

We can now conclude the proof of Proposition \ref{Fubini-like}
by the following lemma.

\begin{lem}
There exists a countable family $(B_{i,j},\chi_{i,j}), (i,j)\in \Nm^2$
of transport beams  such that the images $\chi_{i,j}(B_{i,j})$
cover the set $\mT$.
\end{lem}

\proof
Let $D$ be the closed unit ball in $\Rm^{k-1}$.
Let $\psi_i:D\lto M,i \in \Nm$ be a countable family of 
smooth embeddings such that, for each maximal calibrated curve 
$\gamma:[a,b]\lto M$, the curve 
$\gamma(]a,b[)$ intersects 
the image of $\psi_i$ for some  $i\in \Nm$.
In order to build  such a family of embeddings, let us consider 
a finite atlas $\Theta$ of $M$ composed of charts 
$\theta: B_3\lto M$, where $B_r$ is the open ball of radius
$r$ centered at zero in $\Rm^k$.
We assume that the finite family of open sets $\theta(B_1), \theta \in \Theta$
covers $M$.
For $n=1,\ldots, k$ and $q\in \Qm\cap [-1,1]$, we consider the embedded disk
$\mD_{n,q}\subset B_3$ formed by points $x=(x_1,  \ldots, x_k)\in \bar B_2$
which satisfy $x_n=q$.
The countable family 
$$\theta(\mD_{n,q}); \theta\in \Theta; n=1,\ldots, k;
q\in \Qm\cap [-1,1]
$$
of embedded disks of $M$  
forms a web which intersects all non-trivial curves 
of $M$, hence all transport rays.
We have constructed a countable family of embedded disks which
intersects all transport rays.

For each $(i,j)\in \Nm^2$
let us consider the set $\Omega_{i,j}=D\cap \psi_i^{-1}(\mT_{1/j})$.
Let $a_{i,j}(\omega)$ and  $b_{i,j}(\omega):\Omega_{i,j}\lto \Rm$ be the
functions $-\alpha\circ \psi_i$ and $\beta \circ \psi_i$.
Let $B_{i,j}$ be the set of points $(\omega,s)\subset \Omega _{i,j}\times \Rm$
such that $a_{i,j}(\omega)\leq s\leq b_{i,j}(\omega)$.
To finish, we  define the map $\chi_{i,j}:B_{i,j}\lto M$ by 
$$
\chi_{i,j}(\omega,s)= \gamma_{\psi_i(\omega)}(s).
$$
We claim that, for each $(i,j)\in \Nm^2$, 
the map $\chi_{i,j}$ is Lipschitz, so that
the pair $(B_{i,j},\chi_{i,j})$
is a transport  beam.
In order to prove this claim,
remember that there exists a vector field $E$ on $TM$, the Euler-Lagrange
vector field,  such that the extremals are the projections
of the integral curves of $E$. Because of energy conservation,
this vector field generates a complete flow, denoted by $f_s:TM\rightarrow TM$
for $s\in\RR$. From the fact that the Hamiltonian $H$ is $C^2$ and $f_s$
is in Legendre duality with the Hamiltonian flow, we deduce that
$(s,x,v)\rightarrow f_s(x,v)$ is $C^1$. 
We have 
$$\chi_{i,j}(x,s)=P_M\circ f_s(x,\dot\gamma_{\psi_i(x)}(0))~~
\forall (x,s)\in B_{i,j},$$
where $P_M:TM\rightarrow M$ is the canonical projection
on $M$.
This map is Lipschitz in view of Proposition \ref{Lip}.
If $R$ is a transport ray, it is clear that $R$ is contained in one of the 
images  $\chi_{i,j}(B_{i,j})$. 
\qed

\section{The distinguished  transport plan}\label{proof}
We shall now prove Theorem 2, and hence Theorem 1. Our approach
is based on  remarks in \cite{AmPr:03} and \cite{AKP:04}, however
it seems new, and is  surprisingly simple. Let $\mu$ be a
transport plan that is optimal for the cost $c$, and, among
these optimal transport plans, minimizes the functional
$\int\sigma d\mu$. The existence of such a plan is
straightforward.

\begin{prop}
There exists a set 
\begin{equation}\label{def: Gamma}
\Gamma\subset \cup_{x\in M}\{x\}\times
R_x^+,
\end{equation}
which  is a countable union of compact sets, 
 such that $\mu(\Gamma)=1$ and which  is monotone in the
following sense: If $(x_i,y_i),i\in \{1,\ldots, k\}$ is a finite
family of points of $\Gamma$ and if $j(i)$ is a permutation such
that $y_{j(i)}\in R^+_{x_i}$ then
$$
\sum_{i=1}^k \sigma(x_i,y_{j(i)}) \geq \sum_{i=1}^k
\sigma(x_i,y_i).
$$
\end{prop}

\proof 
Let us consider 
the cost
 function $\zeta$, where
 $\zeta(x,y):M\times M\lto[0,\infty]$
is  the lower semi-continuous function defined by
$
\zeta(x,y)=\sigma(x,y)$ if $u(y)-u(x)=c(x,y)$ and $\zeta
(x,y)=\infty$ if not. 
Note that $\int \zeta d\mu=\int \sigma d\mu$ is finite.
Theorem 3.2 of \cite{AmPr:03} implies the existence of 
a Borel set $\tilde \Gamma$ on which $\mu$ is concentrated, and
which is monotone.
By interior regularity of 
the Borel measure $\mu$,
there  exists a set $\Gamma \subset \tilde \Gamma$
which is a countable union of compact sets and
on which $\mu$ is concentrated.
Being a subset of the monotone set $\tilde \Gamma$, the set $\Gamma$
is itself monotone.
\qed

\begin{defn}{\em
Let $\Lambda$ be the set of points $x\in M$ such that the set
$$\Gamma_x:= \{y\in M:(x,y)\in \Gamma\}
$$
contains more than one point,
where $\Gamma$ is defined in \eqref{def: Gamma}.
}\end{defn}

\begin{lem}
The set $\Lambda$ is Borel measurable.
\end{lem}
\proof
Let $K^n, n\in \Nm$ be an increasing sequence of compact sets such that
$\Gamma=\cup K^n$.
For each $x\in M$, let  $\delta_n(x)$ be the diameter of the compact set 
$K^n_x$ of points $y\in M$ such that $(x,y)\in K^n$.
It is not hard to see that the function $\delta_n(x)$
is upper semi-continuous, hence Borel measurable.
Since $\Gamma_x=\cup_{n\in \Nm} K^n_x$,
we have $\delta(x)=\sup_n \delta_n(x)$, where $\delta(x)$
is the diameter of $\Gamma_x$.
As a consequence, the function $\delta$ is Borel  measurable,
and the set $\Lambda=\{x\in M, \delta(x)>0\}$ is Borel.
\qed

\begin{prop}
We have $\Lambda \subset \mT$, and the intersection $\Lambda \cap
R$ is at most countable for each transport ray $R$.
\end{prop}

\proof If $x\not \in \mT$, then $R_x^+=\{x\}$ hence
$\Gamma_x\subset \{x\}$, and $x\not \in \Lambda$. Let us now
consider a transport ray $R$  that is the image of a maximal calibrated
curve $\gamma:[\alpha,\beta]\lto M$. Let us denote by
$h:[\alpha,\beta]\lto \Rm$ the function $u\circ \gamma$,
which is strictly increasing (Lemma~\ref{increase}). Note
that $\sigma(\gamma(s),\gamma(t))=(h(t)-h(s))^2$ for $s\leq t$ in
$[\alpha,\beta]$. In view of the monotonicity of $\Gamma$, we have
$$
(h(t)-h(s))^2+(h(t')-h(s'))^2 \leq (h(t)-h(s'))^2+(h(t')-h(s))^2
$$
or equivalently
\begin{equation}\label{h}
(h(t)-h(t'))(h(s)-h(s'))\geq 0
\end{equation}
whenever $(\gamma(s),\gamma(t))\in \Gamma$,
$(\gamma(s'),\gamma(t'))\in \Gamma$, $s'\leq t$, and $s\leq t'$.
Following \cite{GaMc:96} or \cite{AmPr:03}, we observe  that this
implies the property:
\begin{equation}\label{order}
(\gamma(s),\gamma(t))\in \Gamma,\;\; (\gamma(s'),\gamma(t'))\in
\Gamma,\;\; s<s'\;\; \Longrightarrow\;\; t\leq t'.
\end{equation}

This property implies that the set of values of $s$ in
$[\alpha,\beta]$ such that $\Gamma_{\gamma(s)}$ contains more
than one element is at most countable. 
Indeed, for all integers $n\geq 1$, let $S_n$ be the set
of values of $s\in[\alpha,\beta]$ such that there exist
$t_1,t_2\in[\alpha,\beta]$ with $(\gamma(s),\gamma(t_1))\in \Gamma$,
$(\gamma(s),\gamma(t_2))\in \Gamma$ and $t_2-t_1\geq 1/n$.
If $s<s'$ are in $S_n$ and if $t_1,t_2$ and $t_1',t_2'$ are as above
with respect to $s$ and $s'$, then $\alpha\leq t_1\leq t_2-1/n
<t_2\leq t_1'\leq t_2'-1/n<t_2'\leq \beta$ and thus
$\beta-\alpha\geq 2/n$. More generally, if $S_n$ contains at least
$j$ points, then $\beta- \alpha\geq j/n$.
As the interval $[\alpha,\beta]$ is bounded,
the set  $S_n$ is finite for all $n$, which leads to the conclusion.
\qed

Theorem 2 can now be proved in a very standard way. In view of
section \ref{fubini}, the set $\Lambda$ has zero Lebesgue measure
in $M$. The set  $Z= M- \Lambda$ is a Borel set of full Lebesgue measure,
$\mu_0(Z)=1$.
Denoting by $\pi_0:M\times M\lto M$ the projection
on the first factor, we observe that 
 the set $\Gamma_Z=\Gamma\cap\pi_0^{-1}(Z)$ is a Borel graph on
which $\mu$ is concentrated
(because $\mu(\Gamma)=1$).
By the easy Proposition 2.1  of
\cite{Am:00}, we conclude that the plan $\mu$ is induced from a
transport map $F$. We then have
$$\int_{M\times M}cd\mu=\int_M c(x,F(x))d\mu_0(x)=K(\mu_0,\mu_1)
=C(\mu_0,\mu_1),$$ so that the map $F$ is optimal for the cost
$c$. This ends the proof of Theorem 2 and Theorem 1.

\small
\bibliographystyle{amsplain}
\providecommand{\bysame}{\leavevmode\hbox
to3em{\hrulefill}\thinspace}

\end{document}